\documentclass[letterpaper,10pt, journal,twoside]{ieeetran}
\IEEEoverridecommandlockouts
\usepackage{amsmath,amssymb,epsfig,color,subfigure,empheq,graphicx,pgfplots}
\usepackage{enumerate,url,wasysym,epstopdf,balance,enumerate}
\usepackage{algorithm,algpseudocode}

\newcommand \be{\mathbf{e}}

\newcommand \bn{\mathbf{n}}
\newcommand \br{\mathbf{r}}
\newcommand \hbr{\hat{\mathbf{r}}}

\newcommand \bv{\mathbf{v}}
\newcommand \bw{\mathbf{w}}
\newcommand \bx{\mathbf{x}}

\newcommand \bp{\mathbf{p}}
\newcommand \bq{\mathbf{q}}

\newcommand \bPhi{\boldsymbol{\Phi}}

\newcommand \bDelta{\boldsymbol{\Delta}}
\newcommand \bdelta{\boldsymbol{\delta}}

\newcommand \bI{\mathbf{I}}

\newcommand \bN{\mathbf{N}}

\newcommand \bX{\mathbf{X}}

\newcommand \bR{\mathbf{R}}

\newcommand \bTheta{{\mathbf{\Theta}}}

\newcommand \tbv{\tilde{\mathbf{v}}}

\newcommand \tbp{\tilde{\mathbf{p}}}
\newcommand \tbq{\tilde{\mathbf{q}}}

\newcommand \tbV{\tilde{\mathbf{V}}}

\newcommand \mcA{\mathcal{A}}

\newcommand \mcD{\mathcal{D}}
\newcommand \mcG{\mathcal{G}}
\newcommand \mcL{\mathcal{L}}

\newcommand \mcF{\mathcal{F}}
\newcommand \mcI{\mathcal{I}}
\newcommand \mcN{\mathcal{N}}
\newcommand \mcM{\mathcal{M}}

\newcommand \mcP{\mathcal{P}}

\newcommand \mcT{\mathcal{T}}

\newtheorem{proposition}{Proposition}
\newtheorem{lemma}{Lemma}

\newtheorem{definition}{Definition}

\newtheorem{assumption}{Assumption}

\title{Graph Algorithms for Topology Identification\\
using Power Grid Probing}

\author{Guido Cavraro and Vassilis Kekatos,~\IEEEmembership{Senior Member,~IEEE}%
\thanks{The authors are with the Bradley Dept. of ECE, Virginia Tech, Blacksburg, VA 24061, USA. Emails: \tt\small{\{cavraro,kekatos\}@vt.edu}. \normalfont Work partially supported by the NSF-CAREER grant 1751085.}
}

\begin{document}
\maketitle
\thispagestyle{empty}

\begin{abstract}
To perform any meaningful optimization task, power distribution operators need to know the topology and line impedances of their electric networks. Nevertheless, distribution grids currently lack a comprehensive metering infrastructure. Although smart inverters are widely used for control purposes, they have been recently advocated as the means for an active data acquisition paradigm: Reading the voltage deviations induced by intentionally perturbing inverter injections, the system operator can potentially recover the electric grid topology. Adopting inverter probing for feeder processing, a suite of graph-based topology identification algorithms is developed here. If the grid is probed at all leaf nodes but voltage data are metered at all nodes, the entire feeder topology can be successfully recovered. When voltage data are collected only at probing buses, the operator can find a reduced feeder featuring key properties and similarities to the actual feeder. To handle modeling inaccuracies and load non-stationarity, noisy probing data need to be preprocessed. If the suggested guidelines on the magnitude and duration of probing are followed, the recoverability guarantees carry over from the noiseless to the noisy setup with high probability.
\end{abstract}

\begin{IEEEkeywords}
Energy systems; identification; smart grid.
\end{IEEEkeywords}

\allowdisplaybreaks

%%%%%%%%%%%%%%%%%%%%%%%%%%%%%%%%%%%%%%%%%%%%%%%%%%%%%%%%
\section{Introduction}\label{sec:intro}
\IEEEPARstart{P}{ower} distribution grids will be heavily affected by the penetration of distributed energy resources. To comply with network constraints, system operators need to know the topologies of their electric networks. Often utilities have limited information on their primary or secondary networks. Even if they know the line infrastructure and impedances, they may not know which lines are energized. 

This explains the recent interest on feeder topology processing. Several works capitalize on the properties of grid data covariance matrices to reconstruct feeder topologies; see e.g., \cite{Deka1}, \cite{BoSch13}. Graphical models have been used to fit a spanning tree relying on the mutual information of voltage data~\cite{WengLiaoRajagopal17}. Tree recovery methods operating on a bottom-up fashion have been devised in \cite{ParkDekaChertkov}; yet they presume non-metered buses have degree larger than two, fail in buses with constant power factor, and lack practical guidelines for handling noisy setups. All the previous approaches build on second-order statistics of grid data. However, since sample statistics converge to their ensemble counterparts only asymptotically, a large number of grid data is typically needed to attain reasonable performance thus rendering topology estimates obsolete.

When the line infrastructure is known, the problem of finding the energized lines can be cast as a maximum likelihood detection problem in~\cite{CaKe2017}, \cite{Sharon12}. Given power readings at all terminal nodes and selected lines, topology identification has also been posed as a spanning tree recovery exploiting the concept of graph cycles~\cite{sevlian2015distribution}. Line impedances are estimated via a total least-squares fit in~\cite{patopa}. Presuming phasor measurements and sufficient input excitation, a Kron-reduced admittance matrix is recovered via a low rank-plus-sparse decomposition~\cite{Ardakanian17}. 

Rather than passively collecting data, an active data acquisition paradigm has been suggested in \cite{BheKeVe2017}: Inverters are commanded to instantaneously vary their power injections so that the operator can infer non-metered loads by processing the incurred voltage profiles. Perturbing the primary controls of inverters to identify topologies in DC microgrids has also been suggested in~\cite{Scaglione2017}. Line impedances have been estimated by having inverters injecting harmonics in \cite{Ciobotaru07}. Instead of load learning, grid probing has been adopted towards topology inference in~\cite{cake2018inverter}, which analyzes topology recoverability via grid probing and estimates the grid Laplacian via a convex relaxation followed by a heuristic to enforce radiality. 

The current work extends \cite{cake2018inverter} on three fronts. First, it provides a graph algorithm for recovering feeder topologies using the voltage deviations induced at all nodes upon probing a subset of them (Section~\ref{sec:id}). Second, topology recoverability is studied under partially observed voltage deviations, an algorithm is devised, and links between the revealed grid and the actual grid are established (Section~\ref{sec:partial}). Third, noisy data setups are handled by properly modifying the previous schemes and by providing probing guidelines to ensure recoverability with high probability (Section~\ref{sec:noisy}). 

%%%%%%%%%%%%%%%%%% GRID MODELING  %%%%%%%%%%%%%%%%%%%%%
\section{Modeling Preliminaries}\label{sec:model}
%Background concepts from graph theory and an approximate electric model are reviewed first.

%\subsection{Tree Graphs}\label{subsec:trees}
Let $\mathcal{G}=(\mcN,\mcL)$ be an undirected tree graph, where $\mcN$ is the set of nodes and $\mcL$ the set of edges $\mcL:=\{(m,n): m,n \in \mcN\}$. A tree is termed rooted if one of its nodes is designated as the root. This root node will be henceforth indexed by 0. In a tree graph, a \emph{path} is the unique sequence of edges connecting two nodes. The nodes adjacent to the edges forming the path between node $m$ and the root are the \emph{ancestors} of node $m$ and form the set $\mcA_m$. Reversely, if $n\in\mcA_m$, then $m$ is a \emph{descendant} of node $n$. The descendants of node $m$ comprise the set $\mcD_m$. By convention, $m \in \mcA_m$ and $m\in \mcD_m$. If $n \in \mcA_m$ and $(m,n)\in \mathcal E$, node $n$ is the \emph{parent} of $m$. A node without descendants is called a \emph{leaf} or \emph{terminal} node. Leaf nodes are collected in the set $\mcF$, while non-leaf nodes will be termed \emph{internal} nodes. For each node $m$, define its \emph{depth} $d_m:= |\mcA_m|$ as the number of its ancestors. The depth of the entire tree is $d_{\mcG}:=\max_{m\in\mcN} d_m$. If $n\in\mcA_m$ and $d_n=k$, node $n$ is the unique $k$-depth ancestor of node $m$ and will be denoted by $\alpha_{m}^k$ for $k=0,\ldots,d_m$. Let also $\mcT_m^k$ denote the subset of the nodes belonging to the subtree of $\mcG$ rooted at the $k$-depth node $m$ and containing all the descendants of $m$.

Our analysis will be built on the concept of the level sets of a node. The \emph{$k$-th level set} of node $m$ is defined as~\cite{cake2018inverter}
\begin{equation}\label{eq:levelset}
\mcN_m^k := \left\{\begin{array}{ll}
\mcD_{\alpha^k_m} \setminus \mcD_{\alpha^{k+1}_m}&,~ k=0,\ldots,d_m-1\\
\mcD_m&,~ k=d_m
\end{array}\right..
\end{equation}
In essence, the level set $\mcN_m^k$ consists of node $\alpha^k_m$ and all the subtrees rooted at $\alpha^k_m$ excluding the one containing node $m$. Since by definition $\mcN^k_m\subseteq \mcD_{\alpha^k_m}$, the level sets satisfy the ensuing properties that will be needed later.

\begin{lemma}[\cite{cake2018inverter}]\label{le:Nmk}
Let $m$ be a node in a tree graph.
\renewcommand{\labelenumi}{(\roman{enumi})}
\begin{enumerate}
\item The node $\alpha^{k}_m$ is the only node in $\mcN_m^k$ at depth $k$; the remaining nodes in $\mcN_m^k$ are at larger depths; 
%\item If $n\in\mcN_m^k$, then $\alpha^{k}_n=\alpha^{k}_m\in\mcN_m^k$;
\item if $n,s\in\mcN_m^k$, then $\alpha^{k}_n=\alpha^{k}_s=\alpha^{k}_m\in\mcN_m^k$;
\item if $m\in \mcF$, then $\mcN_m^{d_m}=\{m\}$; 
\item if $s\in \mcD_n$ and $n\in\mcN$, then $\mcN^{k}_n=\mcN^{k}_s$ for $k < d_n$; and 
\item if $d_m=k$, then $m \in \mcN_m^k$ and $m \notin \mcN_m^\ell$ for $\ell < k$.
\end{enumerate}
\end{lemma}

A radial single-phase distribution grid having $N+1$ buses can be modeled by a tree graph $\mcG=(\mcN,\mcL)$rooted at the substation. The nodes in $\mcN:=\{0,\ldots,N\}$ denote grid buses, and the edges in $\mcL$ lines. Define $v_n$ as the deviation of the voltage magnitude at node $n$ from the substation voltage, and $p_n+jq_n$ as the power injected through node $n$. The voltage deviations and power injections at all buses in $\mcN\setminus \{0\}$ are stacked in $\bv$, $\bp$, and $\bq$. Let $r_\ell+j x_\ell$ be the impedance of line $\ell$, and collect all impedances in $\br+j\bx$. 
The so termed linearized distribution flow (LDF) model approximates nodal voltage magnitudes as~\cite{BW3,Deka1}
\begin{equation}\label{eq:model}
\bv = \bR\bp + \bX\bq
\end{equation}
where $(\bR,\bX)$ are the inverses of weighted reduced Laplacian matrices of the grid~\cite{CaKe2017}. Let $\br_m$ be the $m$-th column of $\bR$, and $R_{m,n}$ its $(m,n)$-th entry that equals~\cite{Deka1}
\begin{equation}\label{eq:entriesR1}
R_{m,n} = \sum_{\substack{\ell = (c,d) \in \mcL \\ c,d\in \mcA_m \cap \mcA_n}} r_{\ell}.
\end{equation}
The entry $R_{m,n}$ can be equivalently interpreted as the voltage drop between the substation and bus $m$ when a unitary active power is injected as bus $n$ and the remaining buses are unloaded. Leveraging this interpretation, the entries of $\bR$ relate to the levels sets in $\mcG$ as follows.

\begin{lemma}[\cite{cake2018inverter}]\label{le:entriesR2}
Let $m$, $n$, $s$ be nodes in a radial grid.
\renewcommand{\labelenumi}{(\roman{enumi})}
\begin{enumerate}
\item if $m\in\mcF$, then $R_{m,m} > R_{n,m}$ for all $n\neq m$; 
\item $n,s\in\mcN_m^k$ for a $k$ if and only if $R_{n,m} = R_{s,m}$; and
\item if $n\in \mcN_m^{k-1}$, $s\in \mcN_m^{k}$, then $R_{s,m}=R_{n,m} + r_{\alpha^{k-1}_m,\alpha^{k}_m}$.
\end{enumerate}
\end{lemma}

%%%%%%%%%%%%%%%%%%%%%%%%%%%%%%%%%%%%%%%%%%%%%%%%%%
\section{Grid Probing using Smart Inverters}\label{sec:probing}
Solar panels and energy storage units are interfaced to the grid via inverters featuring advanced communication, actuation, and sensing capabilities. An inverter can be commanded to shed solar generation, or change its power factor within msec. The distribution feeder as an electric circuit responds within a second and reaches a different steady-state voltage profile. Upon measuring the bus voltage differences incurred by probing, the feeder topology may be identified. Rather than processing smart meter data on a 15- or 60-min basis, probing actively senses voltages on a per-second basis, over which conventional loads are assumed to be invariant.

The buses hosting controllable inverters are collected in $\mcP \subseteq \mcN$ with $P:=|\mcP|$. Consider the probing action at time $t$. Each bus $m\in\mcP$ perturbs its active injection by $\delta_m(t)$ for one second or so. All inverter perturbations $\{\delta_m(t)\}_{m\in\mcP}$ at time $t$ are stacked in the $P$-length vector $\bdelta(t)$. Based on the model in \eqref{eq:model}, perturbations in active power injections incur voltage differences
\begin{equation}\label{eq:dv}
\tbv(t):=\bv(t)-\bv(t-1) = \bR_\mcP \bdelta(t)
\end{equation}
where $\bR_\mcP\in \mathbb{R}^{N\times C}$ is the submatrix obtained by keeping only the columns of $\bR$ indexed by $\mcP$. 

The grid is perturbed over $T$ probing periods. Stacking the probing actions $\{\bdelta(t)\}_{t=1}^T$ and voltage differences $\{\tbv(t)\}_{t=1}^T$ respectively as columns of matrices $\bDelta$ and $\tbV$ yields
\begin{equation}\label{eq:dV}
\tbV = \bR_{\mcP} \bDelta.
\end{equation}
The data model in \eqref{eq:dv}--\eqref{eq:dV} presumes that injections at non-probing buses remain constant during probing and ignores modeling inaccuracies and measurement noise. The practical setup of noisy data is handled in Section~\ref{sec:noisy}.

Knowing $\bDelta$ and measuring $\tbV$ in \eqref{eq:dV}, the goal is to recover the grid connectivity along with line resistances; line reactances can be found by reactive probing likewise. This task of topology identification can be split into three stages: \emph{s1)} finding $\bR_\mcP$ from \eqref{eq:dV}; \emph{s2)} recovering the level sets for all buses in $\mcP$; and \emph{s3)} finding topology and resistances.

At stage \emph{s1)}, if the probing matrix $\bDelta\in\mathbb{R}^{C\times T}$ is full row-rank, then matrix $\bR_{\mcP}$ can be uniquely recovered as $\bR_{\mcP} = \tbV  \bDelta^+$, where $\bDelta^+$ is the right pseudo-inverse of $\bDelta$. Under this setup, probing for $T=P$ times suffices to find $\bR_{\mcP}$.
 
At stage \emph{s2)}, using Lemma~\ref{le:entriesR2} we can recover the level sets for each bus $m\in\mcP$ as follows:
\begin{enumerate}
\item Append a zero entry at the top of the vector $\br_m$.
\item Group the entries of $\br_m$ to find the level sets of node $m$; see Lemma~\ref{le:entriesR2}-(ii).
\item The number of unique values in the entries of $\br_m$ yields the depth $d_m$.
\item Rank the unique values of $\br_m$ in increasing order, to find the depth of each level set; see Lemma~\ref{le:entriesR2}--(iii).
\end{enumerate}

Given the level sets for all $m\in\mcP$, stage \emph{s3)} recovers the grid topology as detailed next.

%%%%%%%%%%%%%%%%%%%%%%%%%%%%%%%%%%%%%%%%%%%%%%%%%%%%%%%
\section{Topology Recovery}\label{sec:id}
By properly probing the nodes in $\mcP$, the matrix $\bR_\mcP$ can be found at stage \emph{s1)}. Then, the level sets for all buses in $\mcP$ can be recovered at stage \emph{s2)}. Nevertheless, knowing these level sets may not guarantee topology recoverability. Interestingly, probing a radial grid at all leaf nodes has been shown to be sufficient for topology identification~\cite[Th.~1]{cake2018inverter}. To comply with this requirement, the next setup will be henceforth assumed; see also \cite{ParkDekaChertkov}.

\begin{assumption}\label{as:FinC}
All leaf nodes are probed, that is $\mcF\subseteq\mcP$.
\end{assumption}

Albeit Assumption~\ref{as:FinC} ensures topology recoverability, it does not provide a solution for stage \emph{s3)}. We will next devise a \emph{recursive} graph algorithm for grid topology recovery. The input to the recursion is a depth $k$ and a maximal subset of probing nodes $\mcP_n^k$ having the \emph{same} $(k-1)$-depth and $k$-depth ancestors. The $(k-1)$-depth ancestor is known and is denoted by $\alpha_n^{k-1}$. The $k$-depth ancestor is known to exist, is assigned the symbol $n$, yet the value of $n$ is unknown for now. We are also given the level sets $\mcN_{m}^k$ for all $m\in\mcP_{n}^k$. The recursion proceeds in three steps.

The \emph{first step} finds the $k$-depth ancestor $n$ by intersecting the sets $\mcN_{m}^k$ for all $m\in\mcP_{n}^k$. The existence and uniqueness of this intersection are asserted next as shown in the appendix.

\begin{proposition}\label{prop:levelsets}
Consider the subset $\mcP_n^k$ of probing nodes located on the subtree rooted at an unknown $k$-depth node $n\in\mcN$. The node $n$ can be found as the unique intersection
\begin{equation}\label{eq:prop1}
\{n\}= \bigcap_{m \in \mcP_n^k} \mcN^k_m.
\end{equation}
\end{proposition}

At the \emph{second step}, node $n$ is connected to node $\alpha_n^{k-1}$. Since $n=\alpha_m^{k}\in\mcN_{m}^k$ and $\alpha_n^{k-1}=\alpha_m^{k-1}\in\mcN_m^{k-1}$, the resistance of line $(n,\alpha_n^{k-1})$ can be found as [Lemma~\ref{le:entriesR2}-(iii)]
\begin{equation}\label{eq:resistance}
r_{\alpha_n^{k-1},n} = r_{\alpha_m^{k-1},\alpha_m^{k}} = R_{\alpha_m^{k},m} - R_{\alpha_m^{k-1},m}
\end{equation}
for any $m\in\mcP_{n}^k$.

The \emph{third step} partitions $\mcP_n^k\setminus \{n\}$ into subsets of buses sharing the same $(k+1)$-depth ancestor. This can be easily accomplished thanks to the next result.

\begin{proposition}\label{pro:levset&subtree}
For nodes $m$ and $m'$ in a tree graph, it holds that $\alpha_m^{k+1} = \alpha_{m'}^{k+1}$ if and only if $\mcN_m^{k} = \mcN_{m'}^{k}$. 
\end{proposition}

Based on Prop.~\ref{pro:levset&subtree} (shown in the appendix), the set $\mcP_n^k\setminus \{n\}$ can be partitioned by grouping buses with identical $\mcN_m^k$'s. The buses forming one of these partitions $\mcP_s^{k+1}$ have the same $k$-depth and $(k+1)$-depth ancestors. Node $n$ was found to be the $k$-depth ancestor. The $(k+1)$-depth ancestor is known to exist and is assigned the symbol $s$. The value of $s$ is found by invoking the recursion with new inputs the depth $(k+1)$, the set of buses $\mcP_s^{k+1}$ along with their $(k+1)$-depth level sets, and their common $k$-depth ancestor (node $n$).

\begin{algorithm}
\caption{Topology Recovery with Complete Data}
\begin{algorithmic}[1]%for adding line numbers
\Require $\mcN$, $\{\mcN_m^{k}\}_{k=0}^{d_m}$ for all $m\in \mcP$.
\State Run \texttt{Root\&Branch}$(\mcP,\emptyset,0)$.
\Ensure Radial grid and line resistances over $\mcN$.
\end{algorithmic}
\textbf{Function} \texttt{Root\&Branch}$(\mcP_n^k,\alpha_n^{k-1},k)$
\begin{algorithmic}[1] %for adding line numbers
\State Identify the node $n$ serving as the common $k$-depth ancestor for all buses in $\mcP_n^k$ via~\eqref{eq:prop1}.
\If{$k > 0$,}
\State Connect node $n$ to $\alpha_n^{k-1}$ with the resistance of \eqref{eq:resistance}.
\EndIf
\If{$\mcP_n^k\setminus\{n\}\neq \emptyset$,}
	\State Partition $\mcP_n^k\setminus\{n\}$ into groups of buses $\{\mcP_s^{k+1}\}$ having identical $k$-depth level sets.
	\State Run \texttt{Root\&Branch}$(\mcP_s^{k+1},n,k+1)$ for all $s$.
\EndIf
\end{algorithmic}\label{alg:full}
\end{algorithm}

To \emph{initialize} the recursion, set $\mcP_n^0=\mcP$ since every probing bus has the substation as $0$-depth ancestor. At $k=0$, the second step is skipped as the substation does not have any ancestors to connect. The recursion \emph{terminates} when $\mcP_n^k$ is a singleton $\{m\}$. In this case, the first step identifies $m$ as node $n$; the second step links $m$ to its known ancestor $\alpha_m^{k-1}$; and the third step has no partition to accomplish. The recursion is tabulated as Alg.~\ref{alg:full}.

%%%%%%%%%%%%%%%%%%%%%%%%%%%%%%%%%%%%%%%%%%%%%%%%%%%%
\section{Topology Recovery with Partial Data}\label{sec:partial}
Although the scheme described earlier probes the grid only through a subset of buses $\mcP$, voltage responses are collected across all buses. This may be unrealistic in distribution grids with limited real-time metering infrastructure, where the operator reads voltage data only at a subset of buses. To simplify the exposition, the next assumption will be adopted.

\begin{assumption}\label{as:MetC}
Voltage differences are metered only in $\mcP$.
\end{assumption}

Under this assumption, the probing model \eqref{eq:dV} becomes
\begin{equation}\label{eq:dVidealred}
\tbV = \bR_{\mcP \mcP} \bDelta
\end{equation}
where now $\tbV$ is of dimension $P \times T$ and $\bR_{\mcP \mcP}$ is obtained from $\bR$ upon maintaining only the rows and columns in $\mcP$. Similar to \eqref{eq:dV}, $\bR_{\mcP\mcP}$ is identifiable if $\bDelta$ is full row-rank. This is the equivalent of stage \emph{s1)} in Section~\ref{sec:probing} under the partial data setup.

\begin{figure}[t]
\centering
\includegraphics[width=0.46\textwidth]{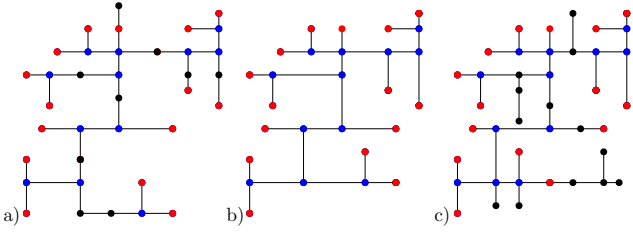}
\caption{\emph{a)} the original IEEE 37-bus feeder; \emph{b)} its reduced equivalent; and \emph{c)} another feeder with the same $\bR_{\mcP\mcP}$. Red nodes are probed; black and blue are not. Blue nodes are internal identifiable nodes comprising $\mcI$.}
\label{fig:redgrid}
\end{figure}

Towards the equivalent of stage \emph{s2)}, since column $\br_m$ is partially observed, only the \emph{metered level sets} of node $m\in\mcP$ defined as $\mcM_m^k := \mcN_m^k \cap \mcP$ can be recovered. The metered level sets for node $m$ can be obtained by grouping the indices associated with the same values in the observed subvector of $\br_m$. Although, the grid topology cannot be fully recovered based on $\mcM_m^k$'s, one can recover a \emph{reduced grid} relying on the concept of internal identifiable nodes; see Fig.~\ref{fig:redgrid}.

\begin{definition}\label{de:I}
The set $\mcI\subset\mcN$ of \emph{internal identifiable} nodes consists of all buses in $\mcG$ having at least two children with each of one of them being the ancestor of a probing bus.
\end{definition}

The reduced grid induced by $\mcP$ can now be defined as the graph $\mcG^r := (\mcN^r, \mcL^r)$ with
\begin{itemize}
\item node set $\mcN^r:=\mcP \cup \mcI$; 
\item $\ell = (m,n) \in \mcL^r$ if $m,n\in\mcN^r$ and all other nodes on the path from $m$ to $n$ in $\mcG$ do not belong to $\mcN^r$; and
\item the resistance of line $\ell = (m,n) \in \mcL^r$ equals the effective resistance between $m$ and $n$ in $\mcG$, that is $r_{mn}^\text{eff}:= (\be_m - \be_n)^\top \bR (\be_m - \be_n)$, where $\be_m$ is the $m$-th canonical vector~\cite{Dorfler13}.
\end{itemize}
In fact, for radial $\mcG$, the resistance $r_{mn}^\text{eff}$ equals the sum of resistances across the $m-n$ path in $\mcG$; see~\cite{Dorfler13}.

Let $\bR^r$ be the inverse reduced Laplacian associated with $\mcG^r$. From the properties of effective resistances, it holds~\cite{Dorfler13}
\begin{equation}\label{eq:sameR}
\bR^r_{\mcP\mcP} = \bR_{\mcP\mcP}.
\end{equation}
In words, the grid $\mcG$ is not the only electric grid having $\bR_{\mcP\mcP}$ as the top-left block of its $\bR$ matrix. The reduced grid $\mcG^r$; the (meshed) Kron reduction of $\mcG$ given $\mcP$; and even grids having nodes additional to $\mcN$ can yield the same $\bR_{\mcP\mcP}$; see Fig.~\ref{fig:redgrid}. Lacking any more detailed information, the grid $\mcG^r$ features desirable properties: i) it connects the actuated and possibly additional identifiable nodes in a radial fashion; ii) it satisfies \eqref{eq:sameR} with the minimal number of nodes; and iii) its resistances correspond to the effective resistances of $\mcG$. Actually, this reduced grid conveys all the information needed to solve an optimal power flow task~\cite{Bolognani2013w}.

The next lemma shows that the number of metered level sets $\mcM_m^k$ coincides with the number of level sets $\mcN_m^k$ for all $m\in\mcP$, so the degrees of probing buses can be reliably recovered even with partial data.

\begin{lemma}\label{lem:metered=level}
Let $\mcG^r = (\mcN^r, \mcL^r)$ be the reduced grid of a radial graph $\mcG$, and let Assumption~\ref{as:FinC} hold true. Then, $\mcN_m^k \cap \mcM_m^k \neq \emptyset$ for all $m \in \mcF$ and $k=1,\ldots,d_m$.
\end{lemma}

\begin{IEEEproof}
Arguing by contradiction, suppose there exists $m\in\mcF$ such that $\mcN_m^k \cap \mcM_m^k = \emptyset$ for some $k\leq d_m$. Since by definition $\alpha_m^k\in\mcN_m^k$, the hypothesis $\mcN_m^k \cap \mcM_m^k = \emptyset$ implies that $\alpha_m^k\notin\mcM_m^k$. Therefore, $\alpha_m^k\notin\mcP$ and the degree of $\alpha_m^k$ is $g_{\alpha_m^k}\geq 3$. The latter implies that $\alpha_m^k$ has at least one child $w\notin\mcA_m$. Let the leaf node $s\in\mcD_w$. Observe that $s$ belongs to both $\mcN^k_m$ and $\mcM^k_m$, contradicting the hypothesis.
\end{IEEEproof}

The next result proved in the appendix guarantees that the topology of $\mcG^r$ is identifiable under Assumption~\ref{as:FinC}.

\begin{proposition}\label{pro:uniqueness_red}
Given a tree $\mcG = (\mcN, \mcL)$ with leaf nodes $\mcF \subseteq \mcN$ and under Assumption~\ref{as:FinC}, its reduced graph $\mcG^r = (\mcN^r, \mcL^r)$ is uniquely characterized by $\{\mcM_m^k \}_{k=0}^{d_m}$ for all $m\in \mcP $, up to different labellings for non-probing nodes.
\end{proposition}

Moving forward to the equivalent of stage \emph{s3) in Section~\ref{sec:probing}}, a three-step recursion operating on metered rather than ordinary level sets is devised next. Suppose we are given the set of probing nodes $\mcP^k_n$ having the same $(k-1)$-depth and $k$-depth ancestors (known and unknown, respectively), along with their $k$-depth metered level sets. 

At the \emph{first step}, if there exists a node $m\in \mcP_n^k$ such that $\mcM_m^k = \mcP^k_n$, then the $k$-depth ancestor $n$ is set as $m$. Otherwise, a non-probing node is added and assigned to be the $k$-depth ancestor. This is justified by the next result.

\begin{proposition}
The root $n$ of subtree $\mcT_n^k$ is a probing node if and only if $\mcM_n^k = \mcP_n^k$.
\end{proposition}

\begin{IEEEproof}%[Proof of Theorem~\ref{th:leaf+id}]
Proving by contradiction, suppose there exists a node $m \in \mcT_n^k$ with $\mcM_m^k = \mcT_n^k \cap \mcP=\mcP_n^k$ and $m\neq n$. Since $m$ is not the root of $\mcT_n^k$, it holds that $d_m > k$, $m \notin \mcM_n^k$, and so $m \notin \mcP_n^k$. If $n$ is a probing node and the root of $\mcT_n^k$, then $d_n=k$ and so $\mcN_n^k=\mcD_n$. Because of this, it follows that $\mcM_n^k = \mcN_n^k \cap \mcP =  \mcD_n \cap \mcP = \mcT_n^k \cap \mcP = \mcP_n^k$.
\end{IEEEproof}

At the \emph{second step}, node $n=\alpha_m^k$ is connected to node $\alpha_n^{k-1}=\alpha_m^{k-1}$. The line resistance can be found through a modified version of \eqref{eq:resistance}. Given any bus $m \in \mcP^k_n$, Lemma~\ref{lem:metered=level} ensures that there exist at least two probing buses $s \in \mcN_m^{k-1}$ and $s' \in \mcN_m^{k}$. Moreover, Lemma~\ref{le:entriesR2}-(ii) guarantees that $R_{\alpha_m^{k-1},m} = R_{s,m}$ and $R_{\alpha_m^{k},m} = R_{s',m}$. Since nodes $m$, $s$, and $s'$ are metered, both $R_{s,m}$ and $R_{s',m}$ can be retrieved from $\bR_{\mcP \mcP}$. Thus, the sought resistance can be found as
\begin{equation}\label{eq:resistance_red}
r_{\alpha_m^{k-1},\alpha_m^k} = R_{\alpha_m^{k},m} - R_{\alpha_m^{k-1},m} =R_{s',m} - R_{s,m}.
\end{equation}

At the \emph{third step}, the set $\mcP_n^k\setminus \{n\}$ is partitioned into subsets of buses having the same $(k+1)$-depth ancestor. This can be accomplished by comparing their $k$-depth metered level sets, as asserted by the next result.

\begin{proposition}\label{pro:levset&subtree_red}
Let $m,m' \in \mcP_n^k$. It holds that $\alpha_m^{k+1} = \alpha_{m'}^{k+1}$ if and only if $\mcM_m^k = \mcM_{m'}^k$.
\end{proposition}

\begin{IEEEproof}
If $\alpha_m^{k+1} = \alpha_{m'}^{k+1}$, then Proposition~\ref{pro:levset&subtree} ensures that $\mcN_m^k = \mcN_{m'}^k$ and so $\mcM_m^k = \mcM_{m'}^k$.

We will show that if $\alpha_m^{k+1} \neq \alpha_{m'}^{k+1}$, then $\mcM_m^k \neq \mcM_{m'}^k$. Since $m,m' \in \mcP_n^k$, it holds that $n = \alpha_m^k = \alpha_{m'}^k$ and
$$\mcM_m^{k} = (\mcD_n \backslash \mcD_{\alpha_m^{k+1}})\cap \mcP,\quad \mcM_{m'}^{k} = (\mcD_n \backslash \mcD_{\alpha_{m'}^{k+1}})\cap \mcP.$$
Because $\mcD_{\alpha_m^{k+1}} \neq \mcD_{\alpha_{m'}^{k+1}}$, $\mcD_{\alpha_m^{k+1}} \cap \mcP \neq 0$, and $\mcD_{\alpha_{m'}^{k+1}} \cap \mcP \neq 0$, it follows that $\mcM_m^{k} \neq \mcM_{m'}^{k}$.
\end{IEEEproof}

The recursion is tabulated as Alg.~\ref{alg:partial}. It is initialized at $k=1$, since the substation is not probed and $\mcM_m^0$ does not exist; and is terminated as in Section~\ref{sec:id}. 

\begin{algorithm}
\caption{Topology Recovery with Partial Data}
\begin{algorithmic}[1]%for adding line numbers
\Require $\mcM$, $\{\mcM_m^{k}\}_{k=1}^{d_m}$ for all $m\in \mcP$.
\State Run \texttt{Root\&Branch-P}$(\mcP,\emptyset,1)$.
\Ensure Reduced grid $\mcG^r$ and resistances over $\mcL^r$.
\end{algorithmic}
\textbf{Function} \texttt{Root\&Branch-P}$(\mcP_n^k,\alpha_n^{k-1},k)$
\begin{algorithmic}[1] %for adding line numbers
\If{ $\exists$ node $n$ such that $\mcM_n^k = \mcP_n^k$,}
\State Set $n$ as the parent node of subtree $\mcT_n^k$.
\Else
\State Add node $n \in \mcI$ and set it as the root of $\mcT_n^k$.
\EndIf
\If{$k>1$,}
\State Connect $n$ to $\alpha_n^{k-1}$ via a line with resistance \eqref{eq:resistance_red}.
\EndIf
\If{$\mcP_n^k\setminus\{n\}\neq \emptyset$,}
	\State Partition $\mcP_n^k\setminus\{n\}$ into groups of buses $\{\mcP_s^{k+1}\}$ having identical $k$-depth metered level sets.
	\State Run \texttt{Root\&Branch-P}$(\mcP_s^{k+1},n,k+1)$ for all $s$.
\EndIf
\end{algorithmic}\label{alg:partial}
\end{algorithm}

%%%%%%%%%%%%%%%%%%%%%%%%%%%%%%%%%%%%%%%%%%%%%%%%%%%%
\section{Topology Recovery with Noisy Data}\label{sec:noisy}
So far, matrices $\bR_\mcP$ and $\bR_{\mcP\mcP}$ have been obtained using the noiseless model of \eqref{eq:dv}. Under a more realistic setup, inverter and voltage perturbations are related as
\begin{equation}\label{eq:dv-noisy}
\tbv(t) = \bR_\mcP \bdelta(t) + \bn(t)
\end{equation}
where $\bn(t)$ captures possible deviations due to non-probing buses, measurement noise, and modeling errors. Stacking $\{\bn(t)\}_{t=1}^T$ as columns of matrix $\bN$, model \eqref{eq:dV} translates to
\begin{equation}\label{eq:dV-noisy}
\tbV = \bR_{\mcP} \bDelta + \bN.
\end{equation}
Under this setup, a least-square estimate can be found as
\begin{equation}\label{eq:LSE2}
\hat\bR_{\mcP} := \arg \min_{\bTheta}\;\|\tbV-\bTheta \bDelta\|_F^2=\tbV\bDelta^+.
\end{equation}

To facilitate its statistical characterization and implementation, a simplified probing protocol is advocated:
\begin{itemize}
\item[\emph{p1})] Every probing bus $m\in\mcP$ perturbs its injection by an identical amount $\delta_m$ over $T_m$ consecutive periods. 
\item[\emph{p2})] During these $T_m$ probing periods, the remaining probing buses do not perturb their injections.
\end{itemize}

Under this protocol, the probing matrix takes the form
\begin{equation}\label{eq:Delta}
\bDelta = \begin{bmatrix}
\delta_1 \be_1 \mathbf{1}_{T_1}^\top & \delta_2 \be_2 \mathbf{1}_{T_2}^\top & \cdots & \delta_P \be_P \mathbf{1}_{T_P}^\top
\end{bmatrix}.
\end{equation}
If at time $t$ node $m$ is probed, the collected $\tbv(t)$ is simply 
\begin{equation}\label{eq:Delta2}
\tbv(t)=\delta_m\br_m+\bn(t).
\end{equation}
Under \eqref{eq:Delta}--\eqref{eq:Delta2}, it is not hard to see that the minimization in \eqref{eq:LSE2} decouples over the columns of $\bR_{\mcP}$. In fact, the $m$-th column of $\bR_{\mcP}$ can be found as the scaled sample mean of voltage differences collected only over the times $\mcT_m:=\left\{\sum_{\tau=1}^{m-1}T_\tau+1,\ldots,\sum_{\tau=1}^{m}T_\tau\right\}$ node $m$ was probed
\begin{equation}\label{eq:sample-mean}
\hbr_m = \frac{1}{\delta_m T_m} \sum_{t\in\mcT_m} \tbv(t).
\end{equation}

To statistically characterize $\hbr_m$, we will next postulate a model for the error term $\bn(t)$ in \eqref{eq:Delta2} as
\begin{equation}\label{eq:noise}
\bn(t) := \bR \tbp(t) + \bX \tbq(t) + \bw(t)
\end{equation}
where $\tbp(t)+j \tbq(t)$ are the injection deviations from non-actuated buses, and $\bw(t)$ captures approximation errors and measurement noise. If $\{\tbp(t),\tbq(t),\bw(t)\}$ are independent zero-mean with respective covariance matrices $\sigma_p^2 \bI$, $\sigma_q^2 \bI$, and $\sigma_w^2 \bI$; the random vector $\bn(t)$ is zero-mean with covariance matrix $\bPhi := \sigma_p^2 \bR^2 + \sigma_q^2 \bX^2 + \sigma_w^2 \bI$.

Invoking the central limit theorem, the estimate $\hbr_m$ can be approximated as zero-mean Gaussian with covariance matrix $\frac{1}{\delta_m^2 T_m}\bPhi$. By increasing $T_m$ and/or $\delta_m$, the estimate $\hbr_m$ can go arbitrarily close to the actual $\br_m$, and this distance can be bounded probabilistically using $\bPhi$. Note however, that $\bPhi$ depends on the unknown $(\bR,\bX)$. To resolve this issue, we resort to an upper bound on $\bPhi$ based on minimal prior information: Suppose the spectral radii $\rho(\bR)$ and $\rho(\bX)$, and the variances $(\sigma_p^2,\sigma_q^2,\sigma_w^2)$ are known; see \cite{Bolognani2013w} for upper bounds. Then, it is not hard to verify that $\rho(\bPhi)\leq \sigma^2$, where $\sigma^2 := \sigma_p^2\rho^2(\bR)+ \sigma_q^2\rho^2(\bX) + \sigma_w^2$. The standard Gaussian concentration inequality bounds the deviation of the $n$-th entry of $\hbr_m$ from its actual value as
\begin{equation}\label{eq:prob}
\Pr \left(|\hat{R}_{n,m} - R_{n,m}| \geq \frac{4 \sigma}{\delta_m \sqrt{T_m}}\right) \leq \pi_0:=6\cdot 10^{-5}.
\end{equation}

Let us now return to stage \emph{s2)} of recovering level sets from the columns of $\bR_\mcP$. In the noiseless case, level sets were identified as the indices of $\br_m$ related to equal values. Almost surely though, there will not be any equal entries in the noisy $\hbr_m$. Instead, the entries of $\hbr_m$ will be concentrated around the actual values. To identify groups of similar values, first sort the entries of $\hbr_m$ in increasing order, and then take the differences of successive sorted entries. A key fact stemming from Lemma~\ref{le:entriesR2}-(iii) guarantees that the minimum difference between the entries of $\br_m$ is larger or equal to the smallest line resistance $r_{\min}$. Hence, if all estimates were confined within $|\hat{R}_{n,m} - {R}_{n,m}| \leq r_{\min}/4$, a difference of sorted $\hat{R}_{n,m}$'s larger than $r_{\min}/2$ would safely pinpoint the boundary between two node groups.

In practice, if the operator knows $r_{\min}$ \emph{a priori} and selects
\begin{equation}\label{eq:Tm}
\delta_m\sqrt{T_m} \geq 16 \sigma/r_{\min}
\end{equation}
the requirement $|\hat{R}_{n,m} - {R}_{n,m}| \leq r_{\min}/4$ will be satisfied with probability higher than $99.95\%$. In such a case and taking the union bound, the probability of correctly recovering all level sets is larger than $1-N^2\pi_0$. The argument carries over to $\bR_{\mcP\mcP}$ under the partial data setup.

%%%%%%%%%%%%%%%%%%%%%% NUMERICAL TESTS  %%%%%%%%%%%%%%%%%%%%%%
\section{Numerical Tests}\label{sec:tests}
Our algorithms were validated on the IEEE 37-bus feeder converted to its single-phase equivalent~\cite{CaKe2017}. Figures~\ref{fig:redgrid}a--\ref{fig:redgrid}b show the actual and reduced topologies that can be recovered under a noiseless setup if all leaf nodes are probed. Setups with complete and partial noisy data were tested. Probing was performed on a per-second basis following the protocol \emph{p1)}--\emph{p2)} of Sec.~\ref{sec:noisy}. Probing buses were equipped with inverters having the same rating as the related load. Loads were generated by adding a zero-mean Gaussian variation to the benchmark data, with standard deviation 0.067 times the average of nominal loads. Voltages were obtained via a power flow solver, and then corrupted by zero-mean Gaussian noise with $3\sigma$ deviation of 0.01\% per unit (pu). Although typical voltage sensors exhibit accuracies in the range of 0.1--0.5\%, here we adopt the high-accuracy specifications of the micro-phasor measurement unit of \cite{mPMU}. 

\begin{table}
\renewcommand{\arraystretch}{1.1}
\caption{Numerical Tests under Full and Partial Noisy Data}
\label{tbl:errors} \centering
\begin{tabular}{|l|l|r|r|r|r|r|}
\hline\hline
& $T_m $ & 1 & $10$ &$20$ & $40$ & $90$ \\
\hline\hline
%IEEE 13-bus & 0.32 & 2.62 & 0.99 & 0.73 & 0.99\\
Alg. 1 & Error Prob. [\%] & 98.5 & 55.3  & 20.9 & 3.1 & 0.2\\
\hline
 & MPE [\%] & 35.1 & 32.5 & 31.2 & 30.9 & 28.5\\
\hline\hline
& $T_m$ & $1$ & $5$ & $10$ & $20$ & $39$\\
\hline\hline
Alg. 2 & Error Prob. [\%] & 97.2 & 45.8 & 26.3 & 18.9 & 0.1\\
\hline
 & MPE [\%] & 18.6 & 16.4 & 15.4 & 14.8 & 13.2\\
\hline
\hline
\end{tabular}
\end{table}

For the 37-bus feeder $r_{\min} = 0.0014$~pu. From the rated $\delta_m$'s; the $r_{\min}$; and \eqref{eq:Tm}, the number of probing actions was set as $T_m = 90$. In the partial data case, the smallest \emph{effective resistance} was $0.0021$~pu, yielding $T_m = 39$. Level sets were obtained using the procedure described in Sec.~\ref{sec:noisy}, and given as inputs to Alg.~\ref{alg:full} and \ref{alg:partial}. The algorithms were tested through 10,000 Monte Carlo tests. Table~\ref{tbl:errors} demonstrates that the probability of error in topology recovery and the mean percentage error (MPE) of line resistances in correctly detected topologies decay gracefully for increasing $T_m$.

\section{Conclusions}\label{sec:conclusions}
To conclude, this letter has put forth an active sensing paradigm for topology identification of inverter-enabled grids. If all lead nodes are probed and voltage responses are metered at all nodes, the grid topology can be unveiled via a recursive algorithm. If voltage responses are metered only at probing buses, a reduced topology can be recovered instead. Guidelines for designing probing actions to cope with noisy data have been tested on a benchmark feeder. Generalizing to multi-phase and meshed grids; coupling (re)active probing strategies; incorporating prior line information; and exploiting voltage phasors are exciting research directions.

%%%%%%%%%%%%%%%%%%%%% APPENDIX %%%%%%%%%%%%%%%%%%%%%%%%%%%
\section*{Appendix}

\emph{Proof of Proposition~\ref{prop:levelsets}.} We will first show that
\begin{equation}\label{eq:lvlset1}
\bigcap_{m \in \mcP_n^k} \mcN^k_m = \bigcap_{m \in \mcT^k_n \cap \mcF} \mcN^k_m.
\end{equation}
By definition $\mcP_n^k=\mcT^k_n \cap \mcP$. Consider a node $w \in \mcP_n^k$ with $w\notin\mcF$. Two cases can be identified. In case i), $w$ equals the root $n$ of subtree $\mcT_n^k$ and hence $\mcN_w^k=\mcD_n=\mcT_n^k$ by the second branch in \eqref{eq:levelset}. Note that $\mcN_s^k \cap \mcD_n = \mcN_s^k$ for all $s \in \mcT_n^k \cap \mcP$. In case ii), node $w$ is different than $n$ and thus Lemma~\ref{le:Nmk}-(iv) implies that $\mcN_w^k = \mcN_s^k$ for all $s \in \mcF \cap \mcD_w$. Either way, it holds
\begin{equation}\label{eq:lvlset2}
\bigcap_{m \in \mcP^k_n} \mcN^k_m = \bigcap_{m \in \left(\mcT^k_n \cap \mcP\right) \setminus \{w\}} \mcN^k_m.
\end{equation}
By recursively applying~\eqref{eq:lvlset2} for each non-leaf probing bus $m$, the equivalence in \eqref{eq:lvlset1} follows. 

From the definition of level sets in \eqref{eq:levelset}, $\mcN_m^k=\mcD_{\alpha^{k}_m}\setminus \mcD_{\alpha^{k+1}_m}$ but $\mcD_{\alpha^{k}_m}=\mcD_n$ since $n$ is the common $k$-depth ancestor for all $m\in \mcP_n^k$. The intersection in the RHS of \eqref{eq:lvlset1} becomes
\begin{equation*}
\bigcap_{m\in \mcT^k_n \cap \mcF} \left(\mcD_n \setminus \mcD_{\alpha^{k+1}_m} \right) = \mcD_n\setminus \bigcup_{m\in \mcT^k_n \cap \mcF}\mcD_{\alpha^{k+1}_m} =\{n\}
\end{equation*}
because $\bigcup_{m\in \mcT^k_n \cap \mcF }\mcD_{\alpha^{k+1}_m}=\mcD_n\setminus \{n\}$.

\emph{Proof of Proposition~\ref{pro:levset&subtree}.}
If $\alpha_m^{k+1} = \alpha_{m'}^{k+1}$, it follows that $\alpha_m^{k} = \alpha_{m'}^{k}$. Then $\mcD_{\alpha^{k+1}_m} = \mcD_{\alpha^{k+1}_{m'}}$ and $\mcD_{\alpha^{k}_m} = \mcD_{\alpha^{k}_{m'}}$. By the definition of the level sets in \eqref{eq:levelset}, it follows that $\mcN_m^{k} = \mcN_{m'}^{k}$.

Conversely, assume that $\mcN_m^{k} = \mcN_{m'}^{k}$. Since $\alpha_m^k$ and $\alpha_{m'}^k$ are the only nodes at depth $k$ respectively in $\mcN_m^{k}$ and $\mcN_{m'}^{k}$ (see Lemma~\ref{le:Nmk}, claim (i)), it follows that $\alpha_m^k = \alpha_{m'}^k$. By the definition of the level sets in \eqref{eq:levelset}, it holds that $\mcN_{m}^{k}=\mcD_{\alpha_{m}^k}\setminus \mcD_{\alpha_m^{k+1}}$, while $\mcD_{\alpha_{m}^{k+1}}\subset \mcD_{\alpha_{m}^{k}}$. Similarly for node $m'$, it holds that $\mcN_{m'}^{k}=\mcD_{\alpha_{m'}^k}\setminus \mcD_{\alpha_{m'}^{k+1}}$ and $\mcD_{\alpha_{m'}^{k+1}}\subset \mcD_{\alpha_{m'}^{k}}$. Since $\mcN_m^{k}=\mcN_{m'}^{k}$ and $\mcD_{\alpha_m^k}=\mcD_{\alpha_{m'}^k}$, it follows that $\mcD_{\alpha_m^{k+1}}=\mcD_{\alpha_{m'}^{k+1}}$ and consequently  $\alpha_m^{k+1} = \alpha_{m'}^{k+1}$.

\emph{Proof of Proposition~\ref{pro:uniqueness_red}}. For the sake of contradiction, assume there exists another reduced grid $\hat \mcG^r = (\mcN^r,\hat \mcL^r)$ with $\hat \mcL^r\neq \mcL^r$ such that $\mcF(\mcG^r) = \mcF(\hat \mcG^r)=\mcF(\mcG)$ and $\{\mcM_w^k(\mcG^r) = \mcM_w^k(\hat \mcG^r)\}_{k=0}^{d_w}$ for all $w \in \mcP$. Note that Lemma~\ref{lem:metered=level} and the latter equality imply that $d_w(\mcG^r)=d_w(\hat{\mcG}^r)$ for all $w\in\mcP$.

Since $\hat \mcG^r \neq \mcG^r$ up to different labelling for non-probing nodes, there exists a subtree $\mcT^{k}_n(\hat \mcG^r)$ with the properties:\\
\emph{1)} It appears both in $\hat \mcG^r$ and $\mcG^r$.\\
\emph{2)} Node $n$ has different parent nodes in $\hat \mcG^r$ and $\mcG^r$, that is $m=\alpha^{k-1}_n(\hat \mcG^r)\neq\alpha^{k-1}_n(\mcG^r)$.\\
\emph{3)} At least one of $\alpha^{k-1}_n(\hat \mcG^r)$ and $\alpha^{k-1}_n(\mcG^r)$ belongs to $\mcP$.

Such a $\mcT^{k}_n(\hat \mcG^r)$ exists and it may be the singleton $\mcT^{k}_n(\hat \mcG^r)=\{n\}$ for $n\in\mcF$. Assume without loss of generality $n\in\mcP$. Based on p3), two cases are identified for $m$.

\emph{Case I}: $m \in \mcP$. We will next show that $d_m(\hat \mcG^r)\neq d_m(\mcG^r)$. From p2) and Lemma~\ref{le:Nmk}-(i), it follows $m\in \mcM_n^{k-1}(\hat \mcG^r)$. On the other hand, property p2) along with the hypothesis $\mcM_s^{k-1}(\hat \mcG^r) = \mcM_s^{k-1}(\mcG^r)$ imply that $d_m(\mcG^r)> k-1=d_m(\hat \mcG^r)$. Lemma~\ref{le:Nmk}-(v) ensures then that $m\in \mcN_m^{k-1}(\hat \mcG^r)$, but $m \notin \mcN_m^{k-1}(\mcG^r)$.

\emph{Case II}: $m\notin \mcP$. Since non-probing buses have degree greater than two in reduced grids, there exists at least one probing node $s$, such that $s\in\mcD_m$, but $s\notin\mcT_n^{k}(\hat\mcG^r)$. Observe that $m,s \in \mcN_n^{k-1}(\hat\mcG^r)$ and $s \in \mcM_n^{k-1}(\hat\mcG^r)$. Let now $w$ be the parent of $n$ in $\mcG^r$. Due to p3), it holds that $w \in \mcP$. Using Lemma~\ref{le:Nmk}-(i), node $w \in \mcM_n^{k-1}(\mcG^r)$ and so $w \in \mcM_n^{k-1}(\hat \mcG^r)$ with possibly $w =s$. Therefore, $d_w(\mcG^r) < d_w(\hat \mcG^r)$ and Lemma~\ref{le:Nmk}-(v) ensures $w \notin \mcN_w^{k-1}(\hat \mcG^r)$ and $w \in \mcN_w^{k-1}(\mcG^r)$.

\bibliographystyle{IEEEtran}
\bibliography{myabrv,power}
\end{document}